\documentstyle[12pt]{article}
\input amssymb.sty
\input epsf

\newtheorem{theorem}{Theorem}[section]
\newtheorem{lemma}[theorem]{Lemma}
\newtheorem{proposition}[theorem]{Proposition}
\newtheorem{corollary}[theorem]{Corollary}
\newtheorem{conjecture}[theorem]{Conjecture}
\newtheorem{definition}[theorem]{Definition}

\begin{document}
\newcommand{\Z}{{\Bbb Z}}
\newcommand{\R}{{\Bbb R}}
\newcommand{\Q}{{\Bbb Q}}
\newcommand{\C}{{\Bbb C}}
\newcommand{\lra}{\longrightarrow}
\newcommand{\lms}{\longmapsto}

\begin{titlepage}
\title{  Multiple polylogarithms, cyclotomy and modular complexes}
\author{A.B. Goncharov }
\date{}
 
\end{titlepage}
\maketitle

\tableofcontents 

\section{Introduction}

 {\bf 1.  Multiple polylogarithms}.   We define them by the  power
series expansion:
\begin{equation} \label{5!5}
Li_{n_{1},...,n_{m}}(x_{1},...,x_{m})  
\quad = \quad 
\sum_{0 < k_{1} < k_{2} < ... < k_{m} } \frac{x_{1}^{k_{1}}x_{2}^{k_{2}}
... x_{m}^{k_{m}}}{k_{1}^{n_{1}}k_{2}^{n_{2}}...k_{m}^{n_{m}}}
\end{equation}
Here $w : = n_1 + ... +n_m$ is called the weight and $m$ the
depth. 

These power series are convergent for $|x_i| <1$, and can 
be continued analytically via the iterated integral presentation given by theorem \ref{BB}. 

The power series (\ref{5!5})  generalize both the 
 classical polylogarithms $Li_n(x)$ (m=1), and      multiple $\zeta$-values  ($x_1 =
... = x_m =1$) :
\begin{equation} \label{lb*}
\zeta(n_{1},...,n_{m}) : = \sum_{0 < k_{1} < k_{2} < ... < k_{m} }
\frac{1}{k_{1}^{n_{1}}k_{2}^{n_{2}}...k_{m}^{n_{m}}} \qquad n_m >1
\end{equation}

The multiple $\zeta$-values were invented and studied by  Euler [E] and 
then   forgotten. They showed up again in such different subjects as quantum groups [Dr] (the Drinfeld associator), Zagier's studies [Z1-2], the Kontsevich integrals for Vassiliev knot invariants, mixed Tate motives over ${\rm Spec} \Z$ [G1-2], and, recently, in computations in quantum field theory [B], [Kr].

The   multiple polylogarithms were studied in   [G1-4]. 
In this paper we  
 investigate them at $N$-th roots of unity: $x_1^N = ... = x_m^N =1$. 
Notice that $Li_1(x) = -\log(1-x)$, so if $\zeta_N$ is a primitive $N$-th root 
of $1$, then $Li_1(\zeta_N )$ is a logarithm of a cyclotomic unit in 
$\Z[\zeta_N, N^{-1}]$. 
In general the supply of numbers we get coincides with the linear
combinations of 
multiple Dirichlet $L$-values
\begin{equation} \label{lb*1}
L(\chi_1,...,\chi_m; n_{1},...,n_{m}) : = \sum_{0 < k_{1} < k_{2} < ... < k_{m} }
\frac{\chi_1(k_1) ...\chi_m(k_m)}{k_{1}^{n_{1}}k_{2}^{n_{2}}...k_{m}^{n_{m}}} 
\end{equation}
They are  periods of  mixed Tate motives over the scheme $S_N:= {\rm Spec} \Z[\zeta_N][\frac{1}{N}]$ 
(see s. 11 of
[G2] and [G4]).

To study these numbers we introduce some tools from
homological algebra (cyclotomic and dihedral Lie algebras,
 modular complex for $GL_m(\Z)$) and geometry (a 
realization of the modular complex in the symmetric space
$SL_m(\R)/SO_m$). To motivate them we start from a conjecture.

{\bf 2. Multiple polylogarithms at roots of unity and the cyclotomic  Lie algebras}.  
Let ${\cal Z}_{\leq w}(N)$ be the $\Q$-vector space spaned by the numbers
\begin{equation} \label{mLv}
\overline Li_{n_1,...,n_m}( \zeta_N^{\alpha_1},...,\zeta_N^{\alpha_m}):= \quad 
 (2 \pi i)^{-w} Li_{n_1,...,n_m}( \zeta_N^{\alpha_1},...,\zeta_N^{\alpha_m})
\end{equation}
Here we may take any branch of $Li_{n_1,...,n_m}(x_1,...,x_m)$. Then 
${\cal Z}(N):= \cup {\cal Z}_{\leq w}(N) $  is  an algebra {\it bifiltered}
by the weight and by the depth.    For example:
\begin{equation} \label{esa}
Li_m(x) \cdot Li_n(y) = \sum_{k_1,k_2 >0}
\frac{x^{k_1}y^{k_2}}{k_1^mk_2^n} = Li_{m,n}(x,y) + 
Li_{m+n}(xy) + Li_{n,m}(y,x)
\end{equation}
(To prove this split the sum over  $k_1 < k_2$, $ k_1=k_2$ and  $ k_1>k_2$).

Denote by  $UC_{\bullet}$  the universal enveloping algebra
of a graded Lie
algebra $C_{\bullet}$. Let $UC_{\bullet}^{\vee}:= \oplus_{n \geq
0}(UC)_{n}^{\vee}$ be its graded dual. It is a commutative Hopf algebra. 

\begin{conjecture} \label{cycle}
There exists a graded Lie algebra $C_{\bullet}(N)$ over $\Q$
   such that one has an isomorphism 
\begin{equation} \label{cycle1}
{\cal Z}(N) \quad = \quad UC_{\bullet}(N)^{\vee}
\end{equation}
of filtered  by the weight on the left and by the degree on the right 
algebras.

b) $H^1_{(n)}(C_{\bullet}(N)) \quad = \quad
K_{2n-1}(\Z[\zeta_N][\frac{1}{N}]) \otimes \Q$. 

c) $C_{\bullet}(1)$ is free graded Lie algebra.
\end{conjecture}

Here $H_{(n)}$ is the degree $n$ part of $H$. 
Notice that $H^1_{(n)}(C_{\bullet}(N))$ is dual to the space of
degree $n$ generators of the Lie algebra $C_{\bullet}(N)$.  

A construction of the Lie algebra $C_{\bullet}(N)$ using the Hodge theory will
be 
outlined 
in s. 3.2. It can be used to deduce conjecture \ref{cycle}   from  some   standard (but extremely dificult!) conjectures in arithmetic
algebraic geometry. 
The simplest abelian quotient $C_{1}(N)$ of  $C_{\bullet}(N)$ is 
the 
group of cyclotomic units in $\Z[\zeta_N, N^{-1}]$, tensored by
$\Q$. We call $C_{\bullet}(N)$ the {\it cyclotomic Lie algebra} of
level $N$,  and  suggest that the "higher cyclotomy theory" should
study its properties.

{\bf Examples}. i) Let $N=1$. Then by the Borel theorem the only
nontrivial 
modulo torsion $K$-groups are 
$K_{4n+1}(\Z)$, which have rank $1$ and correspond to $\overline \zeta(2n+1)$
via the regulator map.

ii) $N=2$: the generators should correspond to $(2\pi i)^{-1}\log 2, \overline \zeta(3), \overline \zeta(5), ... $.

iii) If $N>2$, $n>1$  one has 
$ 
{\rm dim} K_{2n-1}(\Z[\zeta_N, N^{-1}]) \otimes \Q = \frac{\varphi(N)}{2}
$ 
 and the space of generators should correspond to the span over $\Q$ of 
 $ \overline Li_{n }( \zeta_N^{\alpha}
)$.  

{\bf Remark}. Let $\pi^{(l)}_{1}({\Bbb P}^{1} \backslash \{0,1,\infty \})$
 be the $l$-adic completion 
of the fundamental group. One has canonical homomorphism
%$%
\begin{equation} \label{1*q}
\varphi^l: \quad Gal(\overline{\Bbb Q}/ {\Bbb Q}) \longrightarrow Out 
\pi^{(l)}_{1}({\Bbb P}^{1} \backslash \{0,1,\infty \})
%$%
\end{equation}
 It
was studied by P.Deligne, Y.Ihara and others (see [Ih] and references there). 
Conjecture \ref{cycle} for $N=1$ is closely related to some conjectures/questions of P.Deligne [D]
about the image of the map (\ref{1*q})
and V. Drinfeld [Dr] about the structure of the pronilpotent version
of the Grothendieck-Teichmuller group.

The left hand side of (\ref{cycle1}) has an additional structure: the
depth filtration. To study it we proceed as follows.

{\bf 3. The dihedral Lie coalgebra and modular complexes}. 
Let $Z_{\bullet, \bullet}(N)$ be associate graded quotient with respect to
the weight and the depth filtrations of the algebra ${\cal Z}(N)$. 
We reduce it further introducing the bigraded $\Q$-space 
$$
\overline Z_{\bullet, \bullet}(N)  := \quad \frac{Z_{\bullet, \bullet}(N)}
{(Z_{>0, >0}(N) )^2} 
$$
The multiple polylogarithms are multivalued functions, however it is easy 
to show that the projection of (\ref{mLv}) to $\overline Z_{\bullet, \bullet}(N) $   does not depend on the branch we choose.   
So $\overline Z_{w, m}(N)$ is the quotient of the $\Q$-space generated by the
numbers (\ref{mLv}) of weight $w$ and depth $m$ modulo the
subspace generated by the lower weight and depth numbers, and also by the
 products of numbers (\ref{mLv}) of total weight $w$. 

Here is our strategy for investigation of
${\rm dim}_{\Q}\overline Z_{w,m}(N) $. Let $\mu_N$ be the group 
of $N$-th roots of unity. 
In the section 3 a Lie coalgebra ${\cal D}_{\bullet,
\bullet}(\mu_N)$, called the {\it dihedral Lie coalgebra}, is 
explicitely constructed. Namely, the generators of the $\Q$-vector 
space ${\cal D}_{w, m}(\mu_N)$ correspond to   the projections of the numbers 
 (\ref{mLv}) to $\overline  Z_{w, m}(N) $, and the defining relations reflect  
 the known $\Q$-linear relations between these numbers. We prove that 
  $\overline Z_{w,m}(N) $ 
 is a quotient of ${\cal D}_{w, m}(\mu_N)$ (theorem \ref{olmn}). 
 A  really new   data is the 
 cocommutator map $\delta : {\cal D}_{\bullet,
\bullet}(\mu_N)  \lra \Lambda^2 {\cal D}_{\bullet,
\bullet}(\mu_N)$. The dihedral Lie coalgebra is bigraded
 by the weight and the depth. 
 
We want  
to understand the cohomology of the Lie coalgebra ${\cal D}_{\bullet,
\bullet}(\mu_N)$. 
Here is the standard cochain
complex computing  the  cohomology of  ${\cal D}_{\bullet,
\bullet}(\mu_N)$: 
\begin{equation} \label{EEE}
{\cal D}_{\bullet,
\bullet}(\mu_N) \stackrel{\delta }{\lra } \Lambda^2 {\cal D}_{\bullet,
\bullet}(\mu_N)
\stackrel{\delta \wedge id - id \wedge \delta}{\lra } \Lambda^3 {\cal D}_{\bullet,
\bullet}(\mu_N) \stackrel{}{\lra
}  ...
\end{equation}
The first arrow is the cocommutator map, and the others obtained 
via the Leibniz rule. This complex is   bigraded by the weight and
depth. Let $(\Lambda^{\ast}{\cal D}(\mu_N))_{w, m}$ be the subcomplex of
the weight $w$ and depth $m$.
 It is
easy to prove that
\begin{equation} \label{EE}
{\cal D}_{w,1}(\mu_N) \quad = \quad K_{2w-1}(\Z[\zeta_N, N^{-1}]) \otimes \Q \quad \subset 
\quad H^1_{(w)}({\cal D}_{\bullet,\bullet}(N))
\end{equation}

We construct a certain length $m-1$ complex of
$GL_m(\Z)$-modules $M_{(m)}^{\ast}$, 
called the {\it rank $m$ modular
complex}.  For $m=2$ it is identified with 
the chain complex of the  
classical modular triangulation of the hyperbolic plane:

\begin{center}
\hspace{4.0cm}
\epsffile{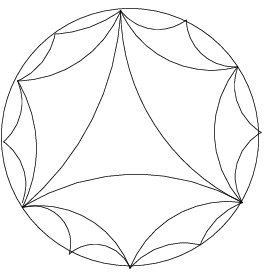}
\end{center}

Let $V_m$ be the standard $m$-dimensional representation of
$GL_m$. Denote by $\Gamma_1(N;m) $ the subgroup of $GL_m(\Z)$
consisting of matrices whose last row is congruent to $(0,...,0,1)$
modulo $N$.

\begin{theorem} \label{D}
a) There exists a   surjective  morphism of complexes
$$
\mu(N)^*_{w,m}: \quad   S^{w-m}V_m \otimes_{\Gamma_1(N; m)}M_{(m)}^{\ast} \quad \lra \quad (\Lambda^{\ast}{\cal D}(\mu_N))_{w,m}
$$

b) It is an
isomorphism if $N=1$, or if $N$ is prime and $w = m$. 
\end{theorem}

A special case of this theorem for $m=2, w=2$ was
proved in [G3]. The case $N=1, m=2$, $w$ is arbitrary was considered 
in [G1] and in s. 4 of [G3]. 

%The map $\mu(N)^*_{w,d}$ is an "almost" isomorphism in general.

We prove in section 6 that the modular complex of rank $3$ is (essentially)
quasiisomorphic to the Voronoi complex of the symmetric
space $SL_3(\R)/SO_3$. Using these results we can estimate 
from above ${\rm dim}\overline Z_{w,m}(N) $ for $m \leq 3$.

{\bf 4. Applications to multiple $\zeta$-values}. Theorem \ref{D}
together with the relation between  the  modular and
Voronoi complexes for $GL_2$ and $GL_3$ mentioned above 
lead to the following result. 
\begin{theorem} \label{C} Let $H_{w,m}$ be the weight $w$, depth $m$ part of 
$H$. Then
 \begin{equation} \label{d1f}
H^i_{(w,2)}({\cal D}_{\bullet,\bullet}(1)) \quad = \quad H^{i-1}(GL_2(\Z), 
S^{w-2}V_2) \quad i=1,2
\end{equation}
  \begin{equation} \label{d2f}
H^i_{(w,3)}({\cal D}_{\bullet,\bullet}(1)) \quad = \quad H^{i}(GL_3(\Z), 
S^{w-3}V_3) \quad i=1,2,3
 \end{equation}
\end{theorem}

This allows us to compute the Euler characterisitc of complexes
$(\Lambda^{\ast}{\cal D}(1))_{w,m}$ for $m=2, 3$. 
Then we find ${\rm
dim}_{\Q}{\cal D}_{w,m}(1)$ for $m=2, 3$ by induction using as a
starting point 
(\ref{EE}), 
which boils 
down to 
 \begin{equation} \label{dep1*}
 {\rm
dim}_{\Q}{\cal D}_{w,1}(1)  = \left\{ \begin{array}{ll}
1 & w :\ {\rm odd} \\ 
 0 & w: \ {\rm even} \end{array} \right.
 \end{equation}

\begin{theorem} \label{46}
 a) If $w$ is odd, then ${\rm dim}\overline Z_{w,2}(1) 
  \quad = \quad {\rm dim}{\cal
D}_{w,2}(1) \quad  = \quad  0$. 

b) If $w$ is even then 
$%\begin{equation} \label{esty1}
{\rm dim}\overline Z_{w,2}(1)  \quad \leq \quad {\rm dim}{\cal
D}_{w,2}(1)\quad  = \quad  [\frac{w-2}{6}].
$%\end{equation}
\end{theorem}
 The part a) goes back to Euler,   b)  was discovered by  Zagier
  [Z2].

 \begin{theorem} \label{B} 
a) If $w$ is even, then ${\rm dim}\overline Z_{w,3}(1) 
 \quad = \quad {\rm dim}{\cal
D}_{w,3}(1) \quad  = \quad  0$. 

b) If $w$ is odd then 
\begin{equation} \label{esty}
{\rm dim}\overline Z_{w,3}(1)  \quad \leq \quad {\rm dim}{\cal
D}_{w,3}(1)\quad  = \quad  [\frac{(w-3)^2-1}{48}]
\end{equation}
\end{theorem} 
Numerical calculations of multiple $\zeta$'s by Zagier
(considerebly extended by Broadhurst in [B]) 
suggested that the estimate (\ref{esty}) is exact. Our results about "motivic"
multiple $\zeta$'s allow to deduce this from standard
conjectures. The details will appear elsewhere.

{\bf Remark}. The philosophy of mixed motives was  the main 
driving force for us. However  constructions and proofs of this paper are 
"elementary", i.e.  do not 
use motives. 
The mixed motives/Hodge structures show up only in s. 3.2 to outline 
the  construction of the cyclotomic Lie algebra, but we do not use s. 3.2
 in the rest of the paper, so the reader may skip it.

{\bf Acknowledgment}. I would like to thank the Max-Planck-Institute (Bonn) 
and University Paris-XI (Orsay) for hospitality and
 support. The  support by the NSF grant DMS-9500010 is gratefully
acknowledged. 

 I  am grateful to J. Bernstein, G. Frey, G. Harder, M. Kontsevich, 
B.B. Venkov, V. Voevodsky, D. Zagier for useful and stimulating discussions.

\section{Properties of multiple polylogarithms}  

{\bf 1. Iterated integral presentation}.  
Set
$$
\int_0^{a_{n+1}} \frac{dt}{a_1-t } \circ ... \circ \frac{dt}{ a_n - t } := 
\int_{0 \leq t_1 \leq ... \leq t_n \leq a_{n+1} } \frac{dt_1}{a_1 - t_1 } 
\wedge ... \wedge\frac{dt_n}{a_n - t_n }
$$
$$
I_{n_1,...,n_m}(a_1:...:a_m:a_{m+1}):=  \int_{0}^{a_{m+1}} \underbrace
 {\frac{dt}{ a_{1} - t} \circ \frac{dt}{t} \circ ... \circ 
\frac{dt}{t}}_{n_{1} \quad  \mbox {times}} \circ   \quad ... \quad \circ 
\underbrace {\frac{dt}{ a_{m} - t} \circ \frac{dt}{t} \circ ... \circ \frac{dt}{t}}_
{n_{m} \quad  \mbox {times}}
$$
 
The following theorem is the key to properties of multiple polylogarithms.
\begin{theorem} \label{BB} 
 $Li_{n_1,...,n_m}(x_1,...,x_m) = I_{n_1,...,n_m}(1: x_1: x_1x_2: ... : x_1...x_m)$.
\end{theorem}
The proof is very easy: 
develope $dt/(a_i - t)$ into a geometric series and integrate. 
If $x_i =1$ we get the Kontsevich formula.
 
{\bf 2. Relations}. {\it The double shuffle relations}.
Set
$$
Li (x_{1},...,x_{m}|t_1,...,t_m ) := \sum_{n_i \geq 1}Li_{n_{1},...,n_{m}}(x_{1},...,x_{m} )t_1^{n_1-1}... t_m^{n_m-1}
$$
$$
I(a_{1}:...:a_{m}:a_{m+1}|t_1,...,t_m ) := \sum_{n_i \geq 1}I_{n_{1},...,n_{m}}(a_{1}:...:a_{m}:a_{m+1})  t_1^{n_1-1}... t_m^{n_m-1}
$$
$$
I^*(a_{1}:...:a_{m}:a_{m+1}|t_1,...,t_m ) := I(a_{1}:...:a_{m}:a_{m+1}|t_1,t_1+t_2,...,t_1+...+t_m )  
$$

Let $\Sigma_{k,n-k}$  be the subset of permutations of $n$ letters $\{1,...,n\}$ consisting of all shuffles of  $\{1,...,k\}$ and $\{k+1,...,n\}$. 
Similar to (\ref{esa}) we see that 
\begin{equation} \label{shuffle1}
Li(x_1,...,x_{k}|t_1,...,t_{k})\cdot
Li(x_{k+1},...,x_{n}|t_{k+1},...,t_{n}) =
\end{equation}
$$
\sum_{ \sigma \in \Sigma_{k,n-k}}Li(x_{\sigma (1)},...,x_{\sigma (n) }|t_{\sigma ( 1) },...,t_{\sigma (n)}) \quad + \quad \mbox{lower depth terms}
$$

\begin{theorem}
\begin{equation} \label{shuffle1}
I^*( a_1:...:a_{k}:1 |t_1,...,t_{k})\cdot
I^*( a_{k+1}:...:a_{n}:1 |t_{k+1},...,t_{n}) =
\end{equation}
$$
\sum_{\sigma \in \Sigma_{k,n-k}}I^*(a_{\sigma (1)},...,a_{\sigma (n)} :1 |t_{\sigma ( 1) },...,t_{\sigma (n)})
$$
 \end{theorem}

{\bf Proof}.  It is not hard to prove the following formula 
\begin{equation} \label{43212}
I^*[a_1:...:a_m:1|t_1,...,t_m] = \int_0^1\frac{ s^{-t_1}}{a_1 -
  s}ds \circ ... \circ \frac{s^{-t_m}}{a_m - s}ds  
\end{equation}
The theorem follows   from this and
the product formula for  iterated integrals.

Here is the simplest case: $
 I_{1 }(x )I_{1 }( y) = I_{1,1}(x,y) + I_{1,1}(y,x )$.
Indeed,
$$
\int_{0}^1 \frac{dt}{t-x} \cdot \int_{0}^1 \frac{dt}{t-y} \quad
= \quad
\int_{0}^1  \frac{dt}{t-x}\circ \frac{dt}{t-y} + \int_{0}^1 \frac{dt}{t-y}\circ
 \frac{dt}{t-x}
$$

For multiple $\zeta$'s these are precisely 
the  relations  
of Zagier, who conjectured that they provide all the relations between the multiple $\zeta$'s 
  
{\it Distribution relations}. From the power series expansion we immediately get
\begin{proposition} \label{2.5}
   If $|x_i| < 1$  and $l$ is a positive integer then  
\begin{equation} \label{5n.new}
 Li (x_{1},...,x_{m}|t_{1},...,t_{m})  \quad  = \quad \sum_{ y^l_i  =  x_i} Li( y_{1},..., y_{m}|lt_{1},...,lt_{m})   
\end{equation}
\end{proposition}

 \section{The dihedral Lie coalgebra of a commutative group} 

{\bf 1. Definitions}. 
Let $G$   be a commutative group. 
We will define a bigraded Lie coalgebra
${\cal D}_{\bullet, \bullet}(G ) = \oplus_{w \geq m \geq 1} {\cal D}_{w, m}(G )$.  
Let us first define  a graded abelian group $\hat {\cal D}_{
\bullet, m}(G  )$. The group $ {\cal D}_{ \bullet, m}(G)$ is its quotient.

Denote by  $C_{m+1}$ the principal homogeneous space of the cyclic
group $\Z/(m+1)\Z$. Let $\Z[C_{m+1}]$ be the  abelian group of
$\Z$-valued functions   on $C_{m+1}$, and   $\Z[C_{m+1}]_0$ is its  quotient by
 constants. Let $Pol_{\bullet}(\Z[C_{m+1}]_0)$ be the
algebra of polynomial functions on $\Z[C_{m+1}]_0$ graded by the
degree. 

Let $D_{m}$ is the dihedral group of symmetries of the $(m+1)$-gon. Set  
\begin{equation}  \label{ddhh}
G^{C_{m+1}}_0:= \quad \{g = (g_1,...,g_{m+1}) \in G^{m+1}| g_1 \cdot ... \cdot g_{m+1} =1 \}
\end{equation} 
We think about the elements of this group as of $m+1$ elements of $G$ 
with product  $1$
located on an oriented circle. 
 Then  the  group $D_{m}$ acts on (\ref{ddhh}). Let $\chi_{m}:D_{m} \to \{\pm 1\}$ 
 be the character trivial on the cyclic 
subgroup and sending the involution to $(-1)^{m+1}$.   Set
$$
\hat {\cal D}_{\bullet +m, m}(G):= \quad \Bigl(\Z[G^{C_{m+1}}_0] \otimes_{\Z} Pol_{\bullet}(\Z[C_{m+1}]_0) \otimes_{\Z} \chi_{m}\Bigr)_{D_{m}} 
$$

Elements of the group  $\hat {\cal D}_{\bullet, m}(G)$ are  presented  by the generating functions
\begin{equation} \label{**}
\{g_0, g_1, ... ,  g_m | t_0:...:t_m\} \quad \mbox{such that } \quad g_0 \cdot  ...  \cdot  g_m = 1 \quad \mbox{and} 
\end{equation}
$$
\{g_0, g_1, ... ,  g_m | t_0:...:t_m\} \quad= \quad\{g_0, g_1, ... ,  g_m | t + t_0:...: t+t_m\}
$$
They satisfy the {\it dihedral symmetry relations}: 
$$
\{g_0,  ... , g_{m-1}, g_m|t_0:t_1:...:t_m\} \quad= \quad\{g_1,  ... , g_m, g_0|t_1:...:t_m: t_0\}  
$$
$$ \{g_0,  ... , g_m|t_0:...:t_m\}\quad = \quad(-1)^{m+1}\{g_m,  ... , g_0|t_m:...:t_0\}
$$
 We picture the elements of $\hat {\cal D}_{\bullet, m}(G)$ as   $m+1$ pairs   $(g_0,t_0), ... , (g_m,t_m)$   located cyclically on an oriented circle:

\begin{center}
\hspace{4.0cm}
\epsffile{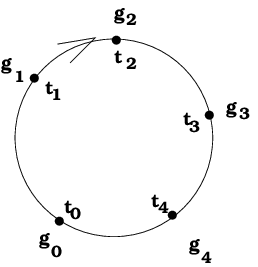}
\end{center}

One has $\hat {\cal D}_{\bullet,m}(G) = \oplus_{w \geq m} \hat {\cal
D}_{w,m}(G)$. We define elements 
\begin{equation} \label{ginf}
\{g_1,...,g_m \}_{n_1,...,n_m}  \in \hat {\cal D}_{w,m}(G), \qquad w =
n_1+...+n_m, \quad n_i \geq 1
\end{equation} 
generating $\hat {\cal D}_{w,m}(G)$ as the coefficients of the generating function:
$$
\{g_0,...,g_m|t_0:....:t_{m}\}:= \sum_{n_i >0} \{g_0, ..., g_{m}  \}_{n_1,...,n_m}
(t_1-t_{0})^{n_1-1}...(t_m -t_{0})^{n_m-1}
$$
Now one can say that $\hat {\cal D}_{w,m}(G)$ is generated by the
elements 
(\ref{ginf}) satisfying the relations imposed by the dyhedral
symmetry. 

We call (\ref{**}) {\it  extended nonhomogeneous dihedral words}  in
 $G$. 
One can  also parametrize the generators    using  the   {\it  extended
 homogeneous dihedral words}:
$$
\{g_0: g_1: ... : g_m|t_0,...,t_m\}, \quad \mbox{where} \quad t_0+ ...+ t_m =0, 
$$  
$$
\{g\cdot g_0:   ... : g\cdot g_m|t_0,...,t_m\} \quad = \quad
\{g_0:  ... : g_m|t_0,...,t_m\}\quad \mbox{ for any } \quad g \in G,
$$
and  the dihedral symmetry holds.  Namely, the duality between the homogeneous and nonhomogeneous extended dihedral words   
is given by  
$$
\{g_0: g_1 : ... :  g_m|t_0,...,t_m\} \longmapsto \{g_0^{-1} g_{ 1}, g_1^{-1} g_{ 2}, ... ,  g_m^{-1} g_{0}|t_0:t_0+t_1:...:t_0+...+t_m \},  
$$
$$
\{g_0, g_1, ... ,  g_m|t_0:...:t_m\} \longmapsto \{g_0: g_0  g_1 : ... :  g_0   ...  g_m|t_1-t_0,t_2-t_1,...,t_0-t_m\}
$$
   
\begin{definition}
${\cal D}_{\bullet, m}(G)$ is the quotient of $\hat {\cal D}_{\bullet, m}(G)$ by  
the following relations: 

a) The double shuffle relations  $(k+l = m, k\geq 1, l \geq 1)$:
$$
 \sum_{\sigma \in \Sigma_{k,l}}\{g_{0}:g_{\sigma(1)} : ... :  g_{\sigma(m)}| t_{0}, t_{\sigma(1)}, ...,  t_{\sigma(m)}\} =0,
$$
$$
 \sum_{\sigma \in \Sigma_{k,l}}
\{x_0, x_{\sigma(1)}, ... ,  x_{\sigma(m)}|t_0: t_{\sigma(1)}: ...:  t_{\sigma(m)} \} =0
$$

b) The distribution relations ($l \in \Z$ and   $|l|$ divides  $|G|$
if the group is finite.)
$$
\{x_0^l, x_{ 1}^l, ... ,  x_{ m}^l| t_0: t_{ 1 }: ...:  t_{ m } \} -
\sum_{y_i^l = x^l_i}\{y_0, y_{ 1}, ... ,  y_{m} | l\cdot t_0: ...:  l\cdot t_{ m} \} =0
$$
 except 
 the relation $\{1\}_1 = \sum_{y^l=1} \{y\}_1$, which is not supposed to  hold.
\end{definition}

{\bf Example}. The distribution relations for $l=-1$ are
 \begin{equation} \label{inversion}
\{x_0^{-1}, x_{ 1}^{-1}, ... ,  x_{ m}^{-1}| t_0: t_{1}: ...:  t_{m} \} =
 \{x_0, x_{ 1}, ... ,  x_{m} | - t_0: ...:  - t_{ m} \}
 \end{equation}

{\bf Remark}. The dihedral symmetry and (\ref{inversion}) follow  
from the double shuffle relations, just copy the proof of theorem \ref{d3} below.

Let us define a cobracket 
$ 
\delta:   \hat {\cal D}_{\bullet, \bullet}(G)  \longrightarrow
 \hat {\cal D}_{\bullet, \bullet}(G)\wedge \hat {\cal
D}_{\bullet, \bullet}(G)
$  
by setting
$$
\delta \{g_0,  ... ,  g_m| t_0:...:t_m\}:= 
$$
$$
\sum_{i=1}^{m-1} \sum_{j=0}^{m} \{g_{j+i+1 }, ... ,g_{j+m},y_{ij}|t_{j+i+1}: ... : t_{j+m+1}    \} \wedge \{x_{ij}, g_{j+1},  ... ,  g_{j+i}|t_j: ... : t_{j+i}  \}  
$$
where indices are modulo $m+1$ and $x_{ij}g_{j+1}  ...   g_{j+i} = 1$, $y_{ij} g_{j+i+1 }  ...  g_{j+m}= 1$. Each term of the formula corresponds to the following procedure: we choose 
an arc on the circle between the two neighboring  distinguished points, and in addition choose 
a distinguished point different from the ends of the arc. 
Then we cut the circle 
in the choosen arc and in the choosen point,  make two naturally oriented circles out of it, and then make  the extended dihedral word on each of the circles out of the initial word in a natural way.

\begin{center}
\hspace{4.0cm}
\epsffile{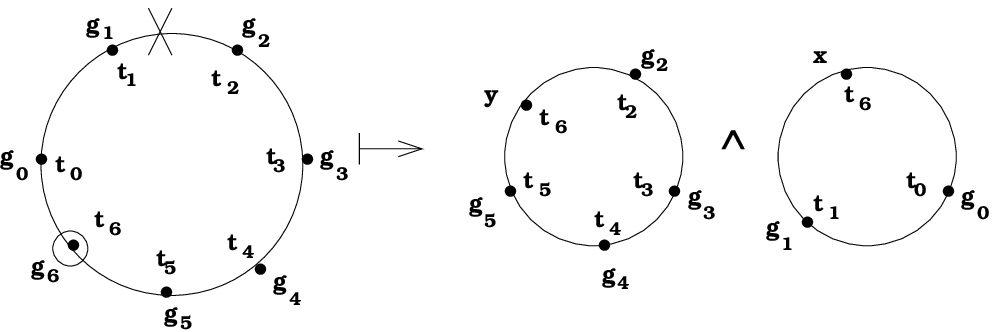}
\end{center}

There is a similar formula for the homogeneous dihedral words,  just exchange $g$'s and $t$'s on the circle.  For example  
$$
\delta \{g_0: g_1 :   g_2 | t_0, t_1, t_2\} \quad = \quad -\{g_2:
g_0 | t_0, -t_0\} \wedge \{ g_1 : g_2 | t_1, -t_1\} - 
$$
$$
-\{g_0: g_1 | t_1,
-t_1\} \wedge \{ g_2 : g_0| t_2, -t_2\} - \{g_1: g_2 | t_2, -t_2\}
\wedge \{ g_0 : g_1 | t_0, -t_0\}
$$

  \begin{theorem}
   $\delta$ provides the structure of   bigraded Lie coalgebra on both 
$\hat {\cal D}_{\bullet, \bullet}(G)$ and ${\cal D}_{\bullet, \bullet}(G)$.
\end{theorem}

{\bf 2.  The cyclotomic Lie algebra and mixed Tate motives over $S_N$}. 
Recall that a mixed Hodge structure is called
a Hodge-Tate structure if all the Hodge numbers $h^{p,q}$ with $p \not
= q$ are zero. The category of mixed $\Q$-Hodge-Tate structures is 
canonically equivalent to the category of finite dimensional comodules 
over a   certain graded pro-Lie  coalgebra  
 ${\cal L}^{HT}_{\bullet}$ over $\Q$ (see   [BGSV], 
[G4]). 
One can attach to the iterated integral related to
$Li_{n_1,...,n_m}(x_1,...,x_m)$ by theorem  \ref{BB}  above an element
$Li^{{\cal H}}_{n_1,...,n_m}(x_1,...,x_m) \in {\cal L}^{HT}_{w}$, the
motivic multiple polylogarithm.  See s.9, 11 of [G2] (or [G1]) where  
${\cal L}^{HT}_{\bullet}$ and a  $w$-framed mixed Hodge-Tate structure 
$Li^{{\cal H}}_{n_1,...,n_m}(x_1,...,x_m)$ were defined. The framed mixed 
Hodge structures corresponding to the  divergent numbers (\ref{mLv}) are 
defined as the limiting Hodge structures, (compare with [G1]  and s.7). For instance   
$Li^{{\cal H}}_1(1) =0$.

\begin{definition} ${\cal C}_{w}(N)$ is the $\Q$-subspace of ${\cal L}^{HT}_{w}$ 
generated by 
the motivic multiple polylogarithms at $N$-th roots of unity of weight $w$.
\end{definition}

\begin{theorem} ${\cal C}_{\bullet}(N) := \oplus_{w \geq 1}{\cal C}_{w}(N)$ is a 
Lie subcoalgebra in ${\cal L}^{HT}_{\bullet}$.
\end{theorem}
We define the cyclotomic Lie algebra  $C_{\bullet}(N)$ as the (graded) dual
of  ${\cal C}_{\bullet}(N)$.

   The hypothetical abelian category of mixed Tate motives over   the scheme $S_N$ is supposed to be canonically equivalent to the category of finite dimensional modules over a graded Lie algebra $L(S_N)_{\bullet}$, called the motivic Lie algebra of 
that scheme, (see [G4] for details). 
$L(S_N)_{\bullet}$ is isomorphic to a free graded Lie algebra
generated 
by the $\Q$-spaces $K_{2n-1}(S_N) \otimes \Q$
in  
degree $n$, $n \geq 1$. 

One can prove that the Lie algebra $C_{\bullet}(N)$ has the same generators, so
it is a
quotient of  $L(S_N)_{\bullet}$.  However in general it is smaller then $L(S_N)_{\bullet}$ since the main result of [G3] implies that $C_{\bullet}(N)$
is not  free when $N$ is sufficiently big. For instance, if $N=p$ is a
prime then 
$$
H_{(2)}^2(C_{\bullet}(p), \Q)\quad = \quad H^1(X_1(p),\Q)_+\quad \oplus \quad \Z[\zeta_p]^* \otimes \Q
$$
where $X_1(p)$ is the modular curve and $+$ means the coinvariants of the complex conjugation.
So it is non zero if $p>3$.

{\it The relation between the dihedral and cyclotomic Lie algberas}. 
Let ${\cal C}_{\bullet, \bullet}(N)$ be the  associate graded  
with respect to the depth filtration on the Lie coalgebra ${\cal
C}_{\bullet}(N)$.

 \begin{theorem}
Assume that $x_i^N=1$. Then the map
$$
\{x_1,...,x_m\}_{n_1,...,n_m} \quad \longmapsto \quad   Li^{{\cal H}}_{n_1,...,n_m}(x_1,...,x_m)
$$
provides a surjective morphism of bigraded 
 Lie coalgebras ${\cal D}_{\bullet,\bullet}(\mu_N) 
 \to {\cal C}_{\bullet, \bullet}(N) $ 
 \end{theorem}
 
 The subspace ${\cal D}_{1, 1}(1)$ is generated by   
$\{1\}_1$. Notice that $Li^{{\cal H}}_1(1) =0$. 
 \begin{conjecture} 
 ${\cal D}_{\bullet, \bullet}(1)/{\cal D}_{1, 1}(1) =
{\cal C}_{\bullet, \bullet}(1)$.
 \end{conjecture}
     But already  ${\cal C}_{1, 1}(N)$ 
     is smaller then ${\cal D}_{1, 1}(N)$ if $N$ is not prime ( take $N=25$).

Let  $(l,N) =1$. 
   One has canonical homomorphism
\begin{equation} \label{1*}
\varphi_N^l: \quad Gal(\overline{\Bbb Q}/ {\Bbb Q}) \longrightarrow 
Out \pi^{(l)}_{1}({\Bbb P}^{1} \backslash \{0, \{\zeta_N^{\alpha}\},\infty \})
\end{equation}
The motivic philosophy suggests that the Lie algebra of the Zariski closure of 
the image of $\varphi_N^l$ is isomorphic 
to the Lie algebra $C_{\bullet}(N)^{\vee} \otimes_{\Q} \Q_l$.

\section{The modular complex}

Let $L_m$ be a  lattice of
rank $m$. The  rank
$m$ modular complex   is a complex of left $GL_m(\Z)$-modules in degrees $[1,m]$,  denoted $M_{(m)}^{\bullet}(L_m)$ or simply $M_{(m)}^{\bullet}$:
$$
M_{(m)}^{1} \stackrel{\partial}{\lra} M_{(m)}^{2} \stackrel{\partial}{\lra}
... \stackrel{\partial}{\lra} M_{(m)}^{m}\\
$$

{\bf  1. The groups $M_{(m)}^{l}$}.   
Let $X$ be a set and  $\Z[[X]]$ be the abelian group of infinite  
$\Z$-linear combinations $\sum n_x \{x\}$ of the generators $\{x\}$, where $x\in X$.  
 The set ${\cal P}_m$ of 
basises $(v_1,...,v_m)$ of $L_m$ 
is a principal homogeneous space  of $GL_m(\Z)$. Then 
 $M_{(m)}^{l} = \Z[[{\cal P} _m ]]/R^l_{(m)}$. The "relations" $R^l_{(m)}$ 
are defined below.

   By definition the generators of the group
 $M_{(m)}^{1}$ are
$$
[v_1,...,v_m], \qquad  \mbox{where $(v_1,...,v_m)$ is a basis of $L_m$}
$$
We will  use   other two notation  for them:
$$
<v_0,...,v_{m}>:= \quad [v_1,...,v_{m}]\qquad \mbox{where} \quad v_0+...+v_m = 0, 
\quad \mbox{and} 
\quad
$$
$$ 
[v_1: ... : v_m] := \quad [v_2 - v_1, v_3-v_2, ... ,  v_m - v_{m-1}, -v_m]
$$
Notice that if $v_1,...,v_m$ is a basis of $L_m$ and $v_0+...+v_m = 0$, then
 omitting any vector of $v_0,...,v_m$ we also get a basis. 

The relations are given by  
 {\it the double shuffle relations}: for any  $1 \leq k \leq m-1$
\begin{equation} \label{sshh1}
s (v_1,...,v_k|v_{k+1},...,v_m):= \quad 
\sum_{\sigma \in \Sigma_{k,m-k}} [v_{\sigma(1)},...,v_{\sigma(m)}] \quad = \quad 0
\end{equation}
 \begin{equation} \label{sshh3}
s (v_1:...:v_k|v_{k+1}:...:v_m):= 
\quad \sum_{\sigma \in \Sigma_{k,m-k}} [v_{\sigma(1)}: ... :v_{\sigma(m)}] \quad =
 \quad 0
\end{equation}

\begin{theorem} \label{d3}
   The double shuffle relations imply the following dihedral  symmetry relations:
$[v_1,...,v_m] = [-v_1,...,-v_m]$, and 
$$
<v_0,...,v_{m}> =  <v_1,...,v_{m},v_0>, \qquad
 <v_0,...,v_{m}>  = 
(-1)^{m+1}<v_{m},...,v_0>
$$
\end{theorem}

{\bf Proof}.  Starting from the shuffle relation $s (v_1| v_2:...:v_m)$:
$$
0 = 
[v_1: v_1+v_2: ... : v_1 + ... + v_m] + ... + [v_2: v_2+v_3: ... : v_1 + ... + v_m]  
$$
we rewrite it as
$ 
[ v_2, ...,  v_{m},  v_{0} ] + s(v_1| v_3, ...,  v_{m},  v_0 ) - 
[ v_3, ...,  v_{m},  v_0,  v_1] =0
$,  getting the cyclic symmetry for the generators $[v_1,... ,v_m]$. 
It is easy see that 
$$
\sum_{k=0}^m (-1)^k s (v_k:v_{k-1}:   ...   : v_1 | v_{k+1}:   ...  : v_m) = 0 
$$
which means that
 $%\begin{equation} \label{sshh}
[v_1: v_2:  ... :v_{m-1}:v_m] \quad = \quad (-1)^{m+1}[v_m: v_{m-1}:  ... :v_2:v_1]
$. %\end{equation} 
A similar result is proved for the $[v_1,...,v_m]$ generators. Using this   we get 
$$
[-v_1, -v_2, ..., -v_{m }] = [v_1 + ... + v_m:  ... : v_{m-1} + v_m: v_m] = 
$$
$$
(-1)^{m+1}[v_m: v_{m-1} + v_m: ... :   v_1 + ... + v_m] = 
(-1)^{m+1}[v_{m-1}, ... ,  v_1, v_0] = [v_0, ... , v_{m-1} ]
$$
So we proved the dihedral symmetry  for  the 
 generators $[v_1,...,v_m]$, and thus for the generators $[v_1:...:v_m]$. 

The generators of the group $M_{(m)}^{l}$ are the symbols
\begin{equation} \label{f}
[v_1,...,v_{m_1}]\wedge ... \wedge [v_{m_{l-1}+1},...,v_{m_l}]  \quad 
\end{equation}
where $m= k_1+...+k_l$, $m_i:= k_1+...+k_i$, $\quad (v_1,...,v_m)$ is
a basis of $L_m$, each of the blocks  $[A_i]:= [v_{m_{i}+1},...,v_{m_{i+1}}]$ 
satisfy   the double shuffle relations, and  the blocks
$[A_i]$ anticommute: 
$$
...\wedge  [A_i]\wedge [A_{i+1}] \wedge ... \quad = \quad - ...\wedge  [A_{i+1}]\wedge [A_{i}]  \wedge ... 
$$

{\bf 2. The differential  $\partial$}. We define a differential 
$
\partial: M_{(m)}^{1} \lra M_{(m)}^{2}
$  by  
$$
\partial:\quad <v_1,...,v_{m+1}> \quad \lms \quad -{\rm Cycle}_{m+1}\Bigl(\sum_{k=1}^{m-1} 
[v_1,...,v_k] \wedge [v_{k+1},...,v_m] \Bigr)
$$
where  $
{\rm Cycle}_{m+1} \Bigl(f(v_1,...,v_{m})\Bigr):= \quad \sum_{i=1}^{m+1}f(v_{i},...,v_{i+m})
$, the indices are modulo $m+1$. Then we extend $\partial$ to the complex $M_{(m)}^{\bullet}$ 
using the Leibniz rule:
$$
\partial( [A_1] \wedge [A_2] \wedge ... \wedge [A_n]) := \quad
\sum_{i=1}^n (-1)^i [A_1] \wedge ... \wedge \partial [A_i] \wedge ... [A_n]  
$$

\begin{theorem} \label{16} The differential $\partial$ is a well defined homomorphism of abelian groups. One has $\partial^2=0$.
\end{theorem}

{\bf Remark}. The modular complex is not the standard cochain complex
of any graded Lie coalgebra.

{\bf 3. Modular complexes for $GL_2$ and $GL_3$}.  Here is explicit description of the modular complexes 
$M_{(m)}^{\bullet}$ for $m=1,2,3$.

{\bf 1. $m=1$}. Then $[v] = [-v]$, and  $M_{(1)}^{1} = \Z$.

{\bf 2. $m=2$}. Then the  modular complex is 
$
M_{(2)}^{1} \stackrel{\partial}{\lra} M^{2}_{(2)}$. 
The group $M_{(2)}^{1}$ is generated by $<v_0, v_1,v_2>$ where $ v_0 +
v_1 + v_2 =0$  and  $(v_1,v_2)$ is a basis in $L_2$. 
The differential is:
\begin{equation} \label{34}
\partial: <v_0, v_1,v_2> \quad \lms \quad -[v_1] \wedge [v_2] -  [v_2] \wedge [v_0] -  [v_0] \wedge [v_1] 
\end{equation}

{\bf 3. $m=3$}. The complex looks as follows:
$
M_{(3)}^{1} \stackrel{\partial}{\lra}
M_{(3)}^{2} \stackrel{\partial}{\lra} M_{(3)}^{3}$. 
The differentials are
$$
\partial: \quad <v_0, v_1, v_2, v_3> \quad  \lms \quad 
$$
$$
- [v_1,v_2] \wedge [v_3]  -  [v_2,v_3] \wedge [v_0]  - [v_3,v_0] \wedge [v_1] - [v_0,v_1] \wedge [v_2] 
$$
$$
 -[v_0] \wedge [v_1,v_2]   -  [v_1] \wedge [v_2,v_3] -  [v_2] \wedge [v_3, v_0]- [v_3] \wedge [v_0, v_1] 
$$

$$
\partial: [v_1,v_2]\wedge [v_3] \quad \lms \quad  
-\Bigl(   [v_1]\wedge [v_2] +  [v_2]\wedge 
[-v_1 - v_2]  +   [-v_1 - v_2]\wedge [v_1]\Bigr)\wedge [v_3]
$$

{\bf 4. Modular complexes and modular cohomology}. 
%\begin{definition} \label{35}
We define the  modular complex $M C^{\ast}(\Gamma, V)$  of a subgroup $\Gamma$ of
$GL_m(\Z)$ with coefficients in a right   $GL_m$-module $V$  as follows:
$$
M C^{\ast}(\Gamma, V):=\quad  V \otimes_{\Gamma} M^{\ast}_{(m)}
$$
Its cohomology are called the modular cohomology $M H^{\ast}(\Gamma, V)$ of $\Gamma$ with coefficients in $V$.
%\end{definition}
If $E$ is the  one element group, then 
$
M C^{\ast}(E,\Z ) = M_{(m)}^{\ast}
$.
 
The group $GL_m(\Z)$ acts  from the right on ${ \Z}[\Gamma \backslash GL_m(\Z)]$, and thus 
 on
$
\Z[\Gamma \backslash GL_m(\Z)] \otimes_{\Z} V
$. 
It is easy to see (Shapiro's lemma) that one has
\begin{equation} \label{w3}
M C^{\ast}(\Gamma, V)=\quad \Bigl( \Z[\Gamma \backslash GL_m(\Z)] \otimes_{\Z} V\Bigr) \otimes_{GL_m(\Z)} M^{\ast}_{(m)}
\end{equation}

{\bf 5. A map from the  modular complex  for $\Gamma_1(N; m)$ with coefficients in $\Z[t_1,...,t_m]$ to the  cochain complex of the dihedral Lie coalgebra of $\mu_N$}.
One has
\begin{equation} \label{3636}
\Gamma_1(N; m) \backslash GL_m(\Z) = \{(\alpha_1,...,\alpha_m)| \quad \alpha_i \in \Z/N\Z, \quad {\rm g.c.d.} (\alpha_1,...,\alpha_m,N)=1\}
\end{equation}

Indeed, the group 
$GL_m(\Z)$ acts from the right on $(\Z/N\Z)^m$ and $\Gamma_1(N; m)$ is the stabilizer of the element 
$(0,... ,0,1)$. The $GL_m(\Z)$-orbit of this element is the right hand side of (\ref{3636}). 

Consider the right $GL_m$-module structure on $\Z[t_1,...,t_m]$ given
by 
$
t_i \cdot g := \sum_{j=1}^m(g^{-1})_{ij}t_j
$.   We will construct a canonical  morphism   of  complexes
$$
\mu(N)_{\bullet, m}^{\ast}: MC^{\ast}(\Gamma_1(N; m), \Z[t_1,...,t_m]) \quad \lra \quad \Lambda^{\ast}\Bigl({\cal D}_{\bullet, \bullet }(\mu_N)\Bigr)_{\mbox{depth = m}}
$$
  where the left hand side is graded by (degree of a polynomial  in $t_i$) $+ m$, and the right hand side by the  weight.

Let  $S^{k}[t_1,...,t_m]$ be  the abelian group of degree $k$ polynomials in $t_1,...,t_m$ with integer coefficients. We will use  (\ref{w3}). 
 Let us  define  first  maps
$$%\begin{equation} \label{a22}
\mu_{w, m}^{1}: MC^{1}(\Gamma_1(N; m), S^{w-m}[t_1,...,t_m]) \quad \lra \quad {\cal D}_{m, w}(\mu_N)
$$%\end{equation}
Choose a basis $(v_1,...,v_m)$ in $V_m$. 
Let $\alpha_0 + \alpha_1 + ... + \alpha_m = 0$. Set
$$
\mu(N)_{\bullet, m}^{1}: \quad \sum_{n_i >0}(\alpha_1,...,\alpha_m) \otimes t_1^{n_1-1} ... t_m^{n_m-1}\otimes  [v_1,...,v_m] \quad \lms 
$$
$$
\{\zeta_N^{\alpha_0}, \zeta_N^{\alpha_1}, ..., \zeta_N^{\alpha_m}|0:t_1:...:t_m\}:= \sum_{n_i >0}\{\zeta_N^{\alpha_0}, \zeta_N^{\alpha_1}, ...,   \zeta_N^{\alpha_m}\}_{n_1,...,n_m}t_1^{n_1-1} ... t_m^{n_m-1} 
$$
\begin{lemma} \label{36q}
The map $\mu_{w, m}^{1}\otimes \Q$ is a   surjective  homomorphism .  
\end{lemma}

{\bf Proof}. It is well defined thanks to the definitions of the modular complex
and the dihedral Lie coalgebra. It is surjective because of the distribution relations.  

Define a map 
$$
\mu(N)_{\bullet, m}^{l}: MC^{l}(\Gamma_1(N; m), \Z[t_1,...,t_m]) \quad \lra \quad {\cal D}_{m, \bullet}(\mu_N)
$$
as follows. Choose a primitive $N$-th root of unity $\zeta_N$. Then 
$$
\mu(N)_{w,m}^{l}: \sum_{n_i>0}( \alpha_1 ,..., \alpha_m ) \otimes t_1^{n_1-1} ... t_m^{n_m-1} \otimes [v_1,...,v_{n_1}]\wedge ... \wedge [v_{n_{l-1}+1},...,v_{n_l}] \lms
$$
$$
\{\zeta_N^{\beta_1},\zeta_N^{\alpha_1}, ..., \zeta_N^{\alpha_{n_1}}|0:t_1:...:t_{n_1}\}\wedge ... \wedge \{\zeta_N^{\beta_l},\zeta_N^{\alpha_{n_{l-1}+1}}, ..., \zeta_N^{\alpha_{n_l}}|0:t_{n_{l-1}+1}:...:t_{n_l}\}
$$
where by definition $\beta_i + \alpha_{n_{i-1}+1}+...+ \alpha_{n_i} =
0$ 
for $i=1,...,l$.

%\begin{theorem} \label{36}
%The map $\mu(N)_{\bullet, \bullet}^{\ast}$ is a surjective 
%homomorphism of bigraded complexes.  It is an isomorphism if $N=1$.
%\end{theorem}

{\bf Proof of theorem \ref{D}}. By the very definitions  the map 
$\mu(N)_{\bullet, \bullet}^*$ 
is a well defined morphism of complexes. It is surjective thanks to 
 lemma \ref{36q}.
 If $N=1$  we have the distribution  relations only for $l= -1$, 
 and they follow  from the double shuffle relations by  
  theorem \ref{d3}.  So the map $\mu(1)_{\bullet, \bullet}^{\ast}$ is an isomorphism. 
  
If $N $ is  a prime and $w=m>1$ there is    one additional distribution relation, 
$\{1,...,1\}_{1,...,,1} = \sum _{x_i^p=1}\{x_1,...,x_m\}_{1,...,,1}$.   
But   
 the shuffle relations  give  
$\sum_{\sigma \in S_m}\{x_{\sigma(1)}, ... , 
x_{\sigma(m)}\}_{1,...,1} =0$, which imply the   distribution relation.

\section{The Voronoi complex}

{\bf 1. Voronoi's cell decomposition of $SL_n(\R)/SO_n$}. Let $Q(V_m)$ be the space of   quadratic forms in an $m$-dimensional vector space $V_m$ over $\R$. 
Denote by ${\cal P}(V_m)$ (resp. $\overline {\cal P}(V_m)$) the cone of positive definite (resp. non negative definite) quadratic forms in  $V_m$.  
Then
$$
{\Bbb H}_m \quad : = \quad SL_m(\R)/SO(m) \quad = \quad {\cal P}(V_m)/\R_+^* 
$$
and $\overline {\cal P}(V_m)/\R_+^*$ is its compactification. For example ${\Bbb H}_{2}$ is the hyperbolic plane.

Any vector $f \in V_m^*$ defines a degenerate non negatively definite quadratic form 
$\varphi(f):= (f,x)^2$. Choose a  lattice $L_m \subset V_m^*$. Let $GL(L_m) \subset GL(V_m)$ be the subgroup preserving the lattice $L_m$. 
Take  the convex hull ${\cal C}(L_m)$ of  the vectors $\varphi(l)$ in 
the cone $\overline {\cal P}(V_m)$ when $l$ runs through all non zero primitive 
vectors of the lattice $L_m$. It is of codimension $1$ in $Q(V_m)$ and has a structure of 
an infinite polyhedra. Its faces are certain convex polyhedras with vertices $\varphi(l_1), ... , \varphi(l_n)$, $l_i \in L_m$.   
 Projecting it  onto ${\cal P}(V_m)/\R_+^*$ we get a
$GL(L_m)$-invariant polyhedral decomposition of   ${\Bbb
H}_{m   }$ called  Voronoi's cell decomposition.    
Set
$$
\varphi(l_1,...,l_n): = \{\lambda_1 \cdot \varphi(l_1) + ... + \lambda_n \cdot 
\varphi(l_n)\}/\R^*, 
\quad \lambda_i \geq 0, \quad \lambda_1 + ... + \lambda_n= 1
$$
The cells of the projection of ${\cal C}(L_m)$ are polyhedras $\varphi(l_1,...,l_n)$ for certain vectors $l_1,...,l_n \in L_m$. They satisfy the condition   ${\rm rk} <l_1,...,l_n> = m$.

The non zero vectors of the lattice $L_m$ minimizing the values of a form $F$ on $V_m^*$ on $L_m \backslash 0$ are called 
{\it the minimal vectors}
of $F$. 
A quadratic form $F$ in $V_m^*$ is called {\it perfect} if the number of minimal vectors of $F$ is at least $\frac{m(m+1)}{2} =   dim Q(V^*_m)$.

Let $s$ be a codimension $1$ face of ${\cal C}(L_m)$. 
Let $h(s)$ be the codimension $1$ subspace in $Q(V_m)$ parallel to the face $s$.

{\bf Voronoi's lemma}. $F \in Q(V^*_m)$ {\it is orthogonal to the subspace $h(s)$ 
if and only if $F$ is a perfect quadratic form. In this case $\{\pm l_1,...,\pm l_n\}$ is the set of minimal vectors for $F$.}

{\bf Proof}. One has $(F, \varphi(l)) = F(l)$. 
Let $(F, x) =c$ be the equation of the hyperplane $h(s)$. Since ${\cal C}(L_m)$ 
is a convex hull it is located in just one of the subspaces $(F, x) < c$ or $(F, x) > c$. Since $(F, \varphi(l)) = F(l)$ could be arbitrary big,   the domain $\{x| (F, x) < c\}$ 
does not intersect ${\cal C}(L_m)$. Further, 
$(F, \varphi(l)) =c$ for any vertex $\varphi(l)$ of the face $\varphi$, so such $l$'s are minimal vectors for $F$. Since the face $\varphi$ is of codimension $1$, 
the number of its vertices is at least $dim Q(V_m)$. So the form $F$ is perfect. The lemma is proved.

Let $v_1,...,v_{m+1}$ vectors of $L_m$ such that $v_1 + ... + v_{m+1} = 0$ 
and $v_1,...,v_m$ is a basis of $L$. 
Set 
\begin{equation} \label{1!}
v_{i,j}:= v_i + v_{i+1} + ... +  v_{j-1}+ v_j, \qquad 1 \leq i, j \leq m+1, \quad i \not = j-1
\end{equation}
 and indices are modulo $m+1$. 
The configuration of vectors $v_{i,j}$ in (\ref{1!})
is linearly equivalent to the configuration of the roots of the root system $A_m$. See fig. 4  for the configuration  of points in $P^2$  
corresponding to the root system $A_2$. 

The convex hull of $\varphi(v_{i,j})$ is a Voronoi cell, called the cell of type $A_m$,  and  
the correseponding perfect form is the quadratic form of 
the root system $A_m$ with the the set of mimimal vectors given by the roots. 

{\bf Voronoi's theorem [V], [M]}. {\it For $m=2,3$ any cell of top dimension in the Voronoi decomposition of ${\Bbb H}_m$ is $GL_m(L_m)$-equivalent to a cell of type $A_m$.}

{\bf 2. The Voronoi complex}. 
Let 
$$
({V}_{\bullet}^{(m)}, d) =({V}_{\bullet}(L_m), d) := \quad {V}^{(m)}_{\frac{m(m+1)}{2} -1} \stackrel{d}{\lra} {V}^{(m)}_{\frac{m(m+1)}{2} -2} 
\stackrel{d}{\lra} ... \stackrel{d}{\lra}  {V}^{(m)}_{m-1}
$$ 
be the complex of (infinite) chains with closed supports  associated with 
the Voronoi decomposition of 
${\Bbb H}_{m}$. We  call it the Voronoi complex of the lattice $L_m$. An isomorphism between lattices lifts to an isomorphism between  the corresponding Voronoi complexes, justifying name {\it the Voronoi complex for $GL_m$}.

\section{Relating the modular and   Voronoi complexes for $GL_2$ and $GL_3$}

The modular complex ${M}_{(m)}^{\bullet}$ is a cohomological complex
placed in degrees $[1,m+1]$. Let us cook up out of him a homological
complex sitting in degrees $[2m, m]$ by setting 
$
{M}^{(m)}_{\bullet}:= {M}_{(m)}^{2m+1 - \bullet}
$. 

{\bf 1. An isomorphism between the modular and Voronoi complexes for  $GL_2$}.  
The Voronoi complex for $GL_2$ looks as follows: 
$
{V}_{\bullet}^{(2)}:= 
{V}^{(2)}_{2} \stackrel{d}{\lra} {V}^{(2)}_{ 1 }
$. 
It is the chain complex of the classical modular triangulation of the
hyperbolic plane, see the figure on page 2

 Define a map of $GL_2(\Z)$-modules  $\psi^{(2)}: M^{(2)}_{\bullet} \lra 
 {V}_{\bullet}^{(2)} $ 
 as follows. Let $v_1,v_2$ be a   basis of $ L_2$,   $
  v_1 + v_2 +v_3  =0$.  Set
$$
[v_1,v_2] \lms \varphi(v_1,v_2,v_3), \qquad 
   [v_1] \wedge [v_2] \lms \varphi(v_1,v_2)
$$

\begin{theorem} \label{thc2} 
  a) The map $\psi^{(2)}$ is an isomorphism of complexes $M^{(2)}_{\bullet} \stackrel{}{\lra} {V}_{\bullet}^{(2)}$. 

 b) Let $\Gamma$ be a  subgroup   of $GL_2(\Z)$. Then for any $GL_2$-module $V$ there are  canonical isomorphisms
$ \quad 
HM_{(2)}^i(\Gamma, V)\otimes {\Q} \quad = \quad H^{i-1}(\Gamma, V)\otimes {\Q}  $.  
\end{theorem}

{\bf Proof}. a) When $(v_1,v_2)$ run  through all basises of the lattice
$L_2$,  the   triangles $\varphi(v_1, v_2, v_3)$  
where $v_1+v_2+v_3
=0$  are cells of type $A_2$, and so by Voronoi's theorem produce 
  all the $2$-cells of  Voronoi's complex for $GL_2$. 
  Since $\psi^{(2)}$ commutes with the differentials, we 
 obviously get an isomorphism of complexes.

b) Voronoi's complex ${V}_{\bullet}^{(2)}$ 
is a resolution of the trivial $GL_2(\Z)$-module $\Z[2]$. 
This resolution is free over   a certain 
finite index subgroup $\Gamma \subset GL_2(\Z)$.

{\bf 2. A quasiisomorphism between the modular and truncated Voronoi complexes for $GL_3$}.
The Voronoi complex for $GL_3$ looks as follows: 
$$
({V}^{(3)}_{\bullet}, d) := \quad
{V}^{(3)}_{5} \stackrel{d}{\lra} {V}^{(3)}_{4} 
\stackrel{d}{\lra}   {V}^{(3)}_{3} \stackrel{d}{\lra}   {V}^{(3)}_{2}
$$ 
  
We will suppose that $v_1,v_2,v_3$ is  a   basis in $L_3$ and 
$$
v_1 + v_2 + v_3 + v_4=0, \qquad v_{12}:= v_1 + v_2, \quad v_{23}:= v_2 + v_3, \quad v_{13}:= v_1 + v_3, \quad ... 
$$

By Voronoi's theorem the $GL_3(\Z)$-orbits of the $5$-symplex $\varphi( v_1,v_2, v_3, v_4,   
v_{12}, v_{23})$ and its faces provide all cells of the Voronoi decomposition for $GL_3$.

Define a map $\psi^{(3)}: M^{(3)}_{\bullet} \lra
{V}_{\bullet}^3/ d {V}_5^{(3)}$ 
 as follows:
$$
[v_1] \wedge [v_2]\wedge [v_3] \lms \varphi(v_1,v_2,v_3)  
$$
$$
[v_1,v_2] \wedge [v_3] \lms \varphi(v_1,v_2,-v_1-v_2, v_3),\quad
$$
\begin{equation} \label{*} 
[v_1,v_2,v_3] \lms \varphi(v_1,v_2, v_3, v_4, v_{12}) - \varphi(v_1,v_2, v_3, v_4,  v_{23})
\end{equation}

\begin{theorem}  \label{thc3} a) The map $\psi^{(3)}$ provides an injective morphism of  complexes 
of $GL_3(\Z)$-modules  
$M^{(3)}_{\bullet} \lra 
{V}_{\bullet}^{(3)}/d {V}_5^{(3)}$. 
It is a  quasiisomorphisms.

b) Let $\Gamma$ be a subgroup   of $GL_3(\Z)$. Then for any $GL_3$-module $V$ there are  canonical isomorphisms
$ 
HM_{(3)}^i(\Gamma, V) \otimes \Q\quad = \quad H^i(\Gamma, V)\otimes \Q, \qquad i \geq 1
$. 
\end{theorem}

{\bf Proof}. b) $=>$ a). 
It is similar to the proof of theorem \ref{thc2} b).

a) Let us show that $\psi^{(3)}$ is a well defined morphism of complexes. 
The map $\psi^{(3)}$ sends the first shuffle relation to zero already in the group ${V}_4^{(3)}$:
$$
\psi^{(3)}: \quad   s (v_1|v_2,v_3) = 
$$
$$
<v_1,v_2,v_3,v_4> + <v_2,v_1,v_3,v_4> + <v_2,v_3,v_1,v_4> \quad \lms 
$$
$$
\varphi(v_1,v_2, v_3, v_4, v_{12}) - \varphi(v_1,v_2, v_3, v_4,  v_{23}) +
\varphi(v_2, v_1,v_3, v_4, v_{12}) - 
$$
$$
- \varphi(v_2,v_1, v_3, v_4,  v_{13}) + \varphi(v_2, v_3,v_1, v_4, v_{23}) - \varphi(v_2, v_3,v_1, v_4,  v_{13}) = 0
$$

The second shuffle relation looks as follows:
$$
s(u_1|u_2:u_3) := [u_1:u_2:u_3] + [u_2:u_1:u_3] +[u_2:u_3:u_1] =
$$
$$
[u_2-u_1,u_3-u_2,-u_3] + [u_1-u_2,u_3-u_1,-u_3] +[u_3-u_2,u_1-u_3,-u_1]
$$
Changing  the variables $v_1:= u_2-u_1, v_2:= u_3-u_2, v_3:= -u_3$ we get 
$$
 <v_1,v_2,v_3,v_4> + <-v_1,v_{12},v_3,v_{41}> + <v_2,-v_{12},-v_4,v_{41}>
$$
The 
maps $\psi^{(3)}$ sends it to the boundary of  Voronoi's  $5$-simplex $\varphi(v_1,v_2, v_3, v_4, v_{12},v_{23})$:
$$
 \varphi(v_1,v_2, v_3, v_4, v_{12}) - \varphi(v_1,v_2, v_3, v_4,  v_{23}) + 
 \varphi(v_1, v_{12}, v_3, v_{23}, v_2) - 
$$
$$
-\varphi(v_1, v_{12}, v_3, v_{23}, v_4) + 
\varphi(v_2,v_{12},v_4,v_{23},v_1) - 
\varphi(v_2,v_{12},v_4,v_{23},v_3) = 
$$
$$
d \varphi(v_1,v_2, v_3, v_4, v_{12},v_{23}) 
$$
The other components of the map $\psi^{(3)}$  are 
obviously group homomorphisms. 
The map $\psi^{(3)}$ respects the differentials,  
and  it is clearly injective. It is an isomorphism in all the degrees except $3$ and $4$. 
 One has
\begin{equation} \label{cocker}
Coker (\psi^{(3)}) = \quad \frac{{  V}_{4}^{(3)}}{\psi^{(3)}(M^{(3)}_{4})} \quad \stackrel{\partial}{\lra}\quad 
\frac{{  V}_{3}^{(3)}}{ \psi^{(3)}(M^{(3)}_{3})}
\end{equation}

Denote by $\{v_1,v_2,v_3,v_4\}$ the set of all unordered $4$-tuple of vectors 
$(v_1,v_2,v_3,v_4)$ in $L_3$  
such that $v_1 + v_2 + v_3 + v_4 = 0$ and 
$(v_1, v_2, v_3)$ is a basis of $L_3$.

\begin{proposition} \label{kker}
One has canonical isomorphisms
\begin{equation} \label{cocker!}
 \frac{V_{4}^{(3)}}{\psi^{(3)}(M^{(3)}_{4})}\quad = \quad \Z[[\{v_1,v_2,v_3,v_4\}]]\quad  = \quad  \frac{V_{3}^{(3)}}{\psi^{(3)}(M^{(3)}_{3})}
\end{equation}
It transforms the differential in (\ref{cocker}) to the identity map on $\Z[[\{v_1,v_2,v_3,v_4\}]]$.  
Therefore the complex $Coker (\psi^{(3)})$ is acyclic. 
\end{proposition}

{\bf Proof}.    
The following  observations about the   Voronoi 
cell decomposition   are  easy to see from  figure 4.

\begin{center}
\hspace{4.0cm}
\epsffile{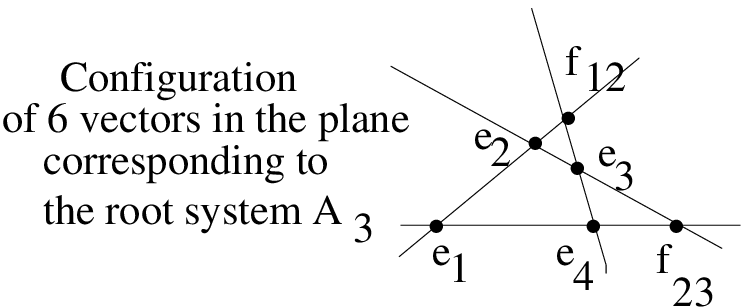}
\end{center}

1) A $3$-dimensional cell of   Voronoi's  decomposition is $GL_3(\Z)$-equivalent to one of the following two:    
a $3$-cell $\varphi(v_1, v_2, -v_{12}, v_3)$, called {\it special $3$-cell}, or
 a $3$-cell $\varphi(v_1, v_2, v_3, v_4)$, called 
{\it generic $3$-cell}.  
 So generic  $3$-cells  are parametrized by the set  $ \{v_1,v_2,v_3,v_4\}$. 

2) A $5$-simplex containing  generic  $3$-cell $\varphi(v_1, v_2, v_3, v_4)$ 
is determined by the dihedral order of the vectors 
$v_1, v_2, v_3, v_4$. So a given generic  $3$-cell   is contained in  three $5$-simplices. 
For the $3$-cell $\varphi(v_1,v_2, v_3, v_4)$ these are 
$$%\begin{equation} \label{three}
\varphi(v_1,v_2, v_3, v_4, v_{12},v_{23}), \quad \varphi(v_2, v_1, v_3, v_4, v_{21},v_{13}), \quad  \varphi(v_2, v_3, v_1, v_4, v_{23},v_{31}) 
$$%\end{equation}

3) The elements (\ref{*})
 are in bijective correspondence with the pairs \linebreak 
$ 
\{ \mbox{generic $3$-cell, a $5$-cell containing it}\}
$.   
The right hand side of (\ref{*}) is the sum of the $4$-cells containing a given generic $3$-cell and contained in a given $5$-cell. 

The first shuffle relation means that the sum of the elements  of type (\ref{*})  corresponding to  $5$-cells
 containing a given generic $3$-cell is zero. 

4) A given $5$-simplex has 
only three generic $3$-cells. 
For the $5$-simplex $\varphi(v_1,v_2, v_3, v_4, v_{12},v_{23})$ they are 
$$
\varphi(v_1, v_2, v_3, v_4), \quad  \varphi(v_{12}, v_4, v_{23}, -v_2), \quad  \varphi(v_{12}, v_3, -v_{23}, -v_1)
$$ 
The sum of the elements (\ref{*}) corresponding to generic 
$3$-cells of a given $5$-simplex is the second shuffle relation. 
It is the boundary of that $5$-simplex. 

The observation  1)  implies the second isomorphism in (\ref{cocker!}). The observations 2) and 3) lead to the first isomorphism in (\ref{cocker!}). 
It is easy to check that the differential is the identity map. The proposition is proved.

\section{Applications to multiple $\zeta$-values}

{\bf 1. Regularization and the map 
${\cal D}_{w,m}(N) \lra \overline Z_{w,m}(N)$}.  
  The iterated integral $I_{n_1,...,n_m}(a_1: ... : a_m:1)$ is divergent 
  if and only if $n_m =1, a_m=1$.     
The power series $Li_{n_1,...,n_m}(x_1,...,x_m)$ are divergent if and only if 
$n_m=1, x_m=1$.

\begin{theorem} \label{olmn}
a)  Let $|x_i | \leq 1 $. Then both the power series
\begin{equation}  \label{ol1}
Li_{n_1,...,n_l,1, ...,1}(x_1,...,x_l, 1- \varepsilon, ... , 1- \varepsilon)
\end{equation} 
   and 
the power series   
$Li_{n_1,...,n_l, 1 , ... , 1}(x_1, ... , x_l,1, ... , 1,1- \varepsilon)
 $ corresponding to  
\begin{equation} \label{ol2}
I_{n_1,...,n_l, 1 , ... , 1}(a_1: ... : a_l:1: ... : 1:1- \varepsilon)
\end{equation}
admit  asymptotic expansions which are polynonials in $  \log
\varepsilon $.  
Their coefficients are explicitely computable sums of multiple polylogarithms.
    
b) The   constant terms of these  
  expansions   differ by    
  lower depth multiple polylogarithms, and the other terms are of
   lower depth. In particular if  $x_i^N=1$  
  they define the same element 
    ${\rm Reg} Li_{n_1,...,n_m}(x_1,...,x_m) \in \overline Z_{\bullet, \bullet}(N)$.

c) The map $\{x_1,...,x_m\}_{n_1,...,n_m} 
\lms {\rm Reg}Li_{n_1,...,n_m}(x_1,...,x_m)$ provides a surjective linear map 
${\cal D}_{w,m}(N) \lra \overline Z_{w,m}(N)$. 
\end{theorem}

{\bf Proof}.    Direct integration gives us     
    $$
     I_{ 1, ...,1}(  1 : ... :1:1- \varepsilon) = 
     \int_{0 < t_1 < ... < t_m < 1 - \varepsilon } \frac{dt_1}{1-t_1}
      \wedge ... \wedge \frac{dt_m}{1-t_m}= 
     \frac{(- \log \varepsilon)^m }{m!} 
    $$  
Computing $  (\sum_{k>0}\frac{(1- \varepsilon)^k}{k})^m  $  
     and using then induction in $m-l $ we get      
    $$
   \frac{(- \log \varepsilon)^m }{m!} \quad = \quad
   Li_{ 1, ...,1}(  1- \varepsilon, ... , 1- \varepsilon)  
        \quad +   
    \sum_{0<i<m}\mbox{(lower depth multiple $\zeta$'s)} (\log \varepsilon)^i
    $$
 Assume $n_l \not = 1$ or $x_l \not = 1$. Applying   the 
 power series product formula to
 $$ 
 Li_{n_1,...,n_l }(x_1,...,x_l) \cdot  
  Li_{ 1, ...,1}(  1- \varepsilon, ... , 1- \varepsilon) 
 $$
 and then to
 $$
 Li_{n_1,...,n_l }(x_1,...,x_l) \cdot  
  Li_{ 1, ...,1}(  1  , ... ,1, 1- \varepsilon),
 $$ 
 and using the induction  on $m-l$, we get parts a) and b) of the theorem.

   c) The power series (resp. the iterated integral) product formulas  clearly hold 
for the asymptotic expansions (\ref{ol1}) (resp. (\ref{ol2})).  
So the shuffle relations are valid. 
The distribution relations for $l>0$ hold for  (\ref{ol1}). 
The distribution relations for $l= -1$ follows from the shuffle relations and 
theorem  \ref{d3}.  
 
  \begin{corollary}
${\rm dim}{\overline Z}(1)_{w,m}  =   0$ if $w+m$ is odd.
\end{corollary}

{\bf Proof}. Indeed, (\ref{inversion}) implies that    ${\rm dim}{\cal D}(1)_{w,m}  =   0$ 
  if $w+m$ is odd.

{\bf 2. Proof of theorem  \ref{C}}. It follows from theorems \ref{D}, \ref{thc2}
 and  \ref{olmn}. 

{\bf 3.  Proof of theorem \ref{46}}.  
 If $w>2$ then  
$H^{0}(GL_2(\Z), S^{w-2}V_2) = 0$, and 
\begin{equation} \label{esty3}
 \sum_{w } {\rm
dim}H^{1}(GL_2(\Z), S^{w-2}V_2) \cdot t^w \quad = \quad \frac{1}{(1-t^4)(1-t^6)}-1
\end{equation}
So by (\ref{d1f}) the generating function for the Euler characteristic of 
the complex 
${\cal D}_{\bullet, 2} \lra \Lambda^2{\cal D}_{\bullet, 1}$ is given by 
(\ref{esty3}). 
 Using    formula (\ref{dep1*}) we get the result. 

{\bf Proof of theorem \ref{B}}. Since $S^{w-3}V_3$ for $w>3$   
is not a self dual $GL_3$-module, the kernel of the restriction of   
$H^i(GL_3(\Z), S^{w-3}V_3)$  to the boundary of the Borel-Serre bordification   
vanishes by a theorem of Borel [BW]. Computing the   boundary contribution 
 to the cohomology 
  we get   
\begin{equation} \label{27}
 H^{i}(GL_3(\Z), 
S^{w-3}V_3)  \quad = \quad \left\{ \begin{array}{ll}
0 &  i=1,2 \\ 
 H_{{\rm cusp}}^{1}(GL_2(\Z),S^{w-2}V_3) &  i=3 \end{array} \right.
\end{equation}
Combining formulas   (\ref{d2f}) and (\ref{27}) we get   
the cohomology of the complex 
${\cal D}_{\bullet, 3} \lra {\cal D}_{\bullet, 2}\otimes {\cal D}_{\bullet, 1}  
\lra \Lambda^3{\cal D}_{\bullet, 1}$.  
Using theorems \ref{46} and \ref{olmn}, 
formula (\ref{dep1*}), and the Euler characteristic argument, we obtain the theorem.

\vskip 3mm \noindent
{\bf REFERENCES}
\begin{itemize}
\item[{[A]}] Ash A.: {\it Cohomolgy of congruence subgroups of
$SL(3,\Z)$}, Math. Ann. 249 (1980) 55-73.
\item[{[BGSV]}] Beilinson A.A., Goncharov A.B., Schechtman V.V., Varchenko A.N.: {\it Aomoto dilogarithms, mixed Hodge structures and motivic cohomolgy}, the Grothendieck Feschtrift, Birkhouser, vol 86, 1990, p. 135-171.
 \item[{[B]}] Borel A., Wallach N.: {\it Continuous cohomology, discrete subgroups and representations of reductive groups} Ann. of Math. Studies, 94, 1980. 
  \item[{[B]}] Broadhurst D.J., {\it On the enumeration of irreducible $k$-fold sums and their role in knot theory and field theory} Preprint hep-th/9604128.
\item[{[E]}] L. Euler: "Opera Omnia," Ser. 1, Vol XV, Teubner,
Berlin 1917, 217-267.
\item[{[G1]}] Goncharov A.B.: {\it Multiple $\zeta$-numbers,
    hyperlogarithms and mixed Tate motives}, Preprint MSRI 058-93, June 1993.
\item[{[G2]}] Goncharov A.B.: {\it Polylogarithms in arithmetic and geometry}, 
Proc. ICM-94, Zurich, 1995, p.374-387.
\item[{[G3]}] Goncharov A.B.: {\it The double logarithm and Manin's
complex for modular curves}. Math. Res. Letters,
 vol. 4. N 5 (1997), pp. 617-636. 
\item[{[G4]}] Goncharov A.B.: {\it Multiple polylogarithms at roots of
unity and motivic Lie algebras},  in Preprint MPI-62/97, Proc. of Arbeitstagung,  Bonn,
June 1997.  
\item[{[G4]}] Goncharov A.B.: {\it Mixed elliptic motives}, Proc. of 60-th Durham symposium ``Galois groups in arithmetic algebraic geometry'', 1998. 
\item[{[D]}] Deligne P.: {\it Le group fondamental de la droite projective moine trois points}, In: Galois groups over $\Bbb Q$. Publ. MSRI, no. 16 (1989) 79-298.  
\item[{[Dr]}] Drinfeld V.G.: {\it On quasi-triangular quasi-Hopf algebras and some group related to  Gal$(\overline{\Bbb Q}/\Bbb Q)$}. Leningrad Math. J., 1991. 
\item[{[Ih]}] Ihara Y.: {\it Braids, Galois groups, and some arithmetic functions}. Proc. Int. Congress of Mathematicians in Kyoto, (1990). 
\item[{[Ih1]}] Ihara Y.: {\it Galois representation arising from
${\Bbb P}^1\backslash \{0,1,\infty\}$ and Tate twists of even degree},
In: Galois groups over $\Bbb Q$. Publ. MSRI, no. 16 (1989).   
 \item[{[Kr]}] Kreimer D.: {\it Renormalisation and knot theory} 
Preprint U. of Mainz, 1996.
\item[{[M]}] Martinet J.: 
{\it Les reseaux parfaits des espaces euclidiens}, (1997), Masson. 
\item[{[V]}] Voronoi G.: {\it Nouvelles applications des
param\'etres continus \'a la th\'eorie des formes quadratiques, I
 %I Sur quelques propri\'et\'es des formes quadratiques positives parfaites
},
J. Reine Angew. Math. 133 (1908), 97-178.
\item[{[Z1]}] Zagier D.: {\it Values of zeta functions and their
    applications}, Proceedings of the First Europian Congress of
  Mathematicians, Paris. 1994, vol 1. 
 \item[{[Z2]}] Zagier D.: {\it Periods of modular forms,
traces of Hecke operators, and multiple zeta values}, RIMS Kokyuroku 843 (199)
162-170.
 \end{itemize}

Dept of Mathematics, Brown University, Providence RI 02912,
USA. e-mail: sasha@math.brown.edu

\end{document}